\documentclass[11pt,a4paper]{article}

\usepackage{lineno,hyperref}
\modulolinenumbers[5]
\usepackage[english]{babel}
\usepackage[a4paper,tmargin=3cm,bmargin=3cm,rmargin=2.5cm,lmargin=2.5cm]{geometry}
\usepackage{amsfonts,amsmath,amssymb}
\usepackage[T1]{fontenc}
\numberwithin{equation}{section}
\newcommand{\disp}{\displaystyle}

\newtheorem{thm}{Theorem}[section]
\newtheorem{cor}[thm]{Corollary}
\newtheorem{lem}[thm]{Lemma}
\newtheorem{prop}[thm]{Proposition}
\newtheorem{defn}[thm]{Definition}

\numberwithin{equation}{section}
\begin{document}
\pagestyle{myheadings}{\markboth{\scriptsize\it $\beta-$Flatness-
Multiplicity Results.} {\scriptsize\it $\beta-$Flatness-
Multiplicity Results.}}
\date{}
\title{\LARGE The $\beta-$Flatness Condition in CR
Spheres\\
Multiplicity Results}

\author{ Najoua Gamara, Boutheina  Hafassa and Akrem Makni }
\maketitle
\begin{center}
\textit{College of Science, Taibah University, KSA}\\
\textit{ University Tunis El Manar, University Campus 2092, Tunisia}
\end{center}
 \date{ }

\begin{abstract}
 We give  multiplicity results for the problem of prescribing
the scalar curvature on Cauchy-Riemann spheres under $\beta-$
flatness condition. To give a lower bound for the number of
solutions, we use Bahri's methods based on the theory of critical
points at infinity and a Poincare'-Hopf type formula.
\end{abstract}
\rule[1mm]{\linewidth}{0.2mm}

\section{Introduction}

\par

In an earlier paper we discussed existence results for the problem
of prescribing the Webster scalar curvature on the
$3$-Cauchy-Riemann sphere, under $\beta$- flatness condition, $2
\leq \beta < 4 $. The purpose of the present paper, is to study
multiplicity results for this problem.

Let $\mathbb{S}^{3}$ be the unit sphere of $\mathbb{C}^{2}$ endowed
with its standard contact form $\theta _{1},$ and
$K:\mathbb{S}^{3}\rightarrow \mathbb{R}$ be a given $C^{2}$ positive
function. The problem of finding a contact form $\theta$ on
$\mathbb{S}^{3}$  conformal to $\theta _{1}$ admitting the function
$K$ as Webster scalar curvature, is equivalent to the resolution of
the following semi-linear equation:
\begin{eqnarray}\label{1}
\left\{
\begin{array}{l}
L_{\theta _{1}}u=K\;u^{3}\quad \text{ on }\mathbb{S}^{3} \\
\quad\;\, u>0
\end{array}
\right.
\end{eqnarray}
where $L_{\theta_{1}}=\,4\Delta_{\theta _{1}} + R_{\theta _{1}}$, is
the conformal laplacian of $\mathbb{S}^{3}$.

%%%%%%%%%%%%%%%%%%%%%%%%%%%%%%%%%%%%%%%%%%%%%%%%%%%%%%%%%%%%%%%%%%%%%%%%%%%%%%%%%%%%%%%%%%%%%%%%%%%%%5

\par Using the CR equivalence $F$
induced by the Cayley Transform (see Definition 2.1 below) between
$\mathbb{S}^{3}$ minus a point and the Heisenberg group
$\mathbb{H}^{1}$, equation (\ref{1}) is equivalent up to an influent
constant to
\begin{eqnarray}\label{2}
\left\{
\begin{array}{l}
4\Delta_{{\mathbb{H}}^{1}}u=\tilde{K}\;u^{3}\quad
\text{ on }\;{\mathbb{H}}^{1}\,\,, \\
\qquad\qquad\quad u>0
\end{array}
\right.
\end{eqnarray}
where $\;\Delta_{\mathbb{ H}^{1}}$ is the sub laplacian of $\mathbb{ H}^{1}$ and $\tilde{K} = K\circ F^{-1}.$

\par In order to give our new multiplicity  results for problem $\eqref{1},$
where the prescribed  function $K$ satisfies a $\beta-$flatness
condition near its critical points. We will use the same techniques
displayed in \cite{Gamara Riahi multiplicity} which are based on an
adaptation to the Cauchy-Riemann settings  of Bahri's work.  These
techniques  were first introduced by Bahri and Coron in \cite{Bahri
Coron}: we have  to study the critical points at infinity of the
associated variational problem, by computing their total Morse
index. Then, we compare this total index
to the Euler characteristic of the space of variation.\\

To state our results, we set up the following conditions and notations.\\
Let $G(a,)$ be a Green's function for $L$ at $a\in\mathbb{S}^{3}$.\\
We denote by \begin{eqnarray*} \mathcal K= \Big\{(\xi_i)_{(1\leq
i\leq r)},\;\text{such that}\; \nabla K(\xi_i)=0\Big\}
\end{eqnarray*}
the set of all  critical points of $K.$ We say that $K$ satisfies
 the \textbf{$\beta-$flatness condition}  if for all $\xi_i
\in \mathcal K,$
 there exist\\
 $$\beta=\beta(\xi_i)\;\;\text{ and} \;\; b_{1}=b_{1}(\xi_i),\ b_{2}=b_{2}(\xi_i),\ b_{0}=b_{0}(\xi_i)\ \in\mathbb R^{\ast}$$
 such that in some pseudo hermitian normal coordinates system centered at  $\xi_i,$ we have
\begin{eqnarray}\label{3'}
K(x)=K(0)+b_1|x_{1}|^{\beta}+b_{2}|x_{2}|^{\beta}+
b_{0}|t|^{\frac{\beta}{2}}+\mathcal{R}(x).
\end{eqnarray}
Where $\disp\sum_{k=1}^{2}b_{k}+\kappa b_{0}\neq 0,
\disp\sum_{k=1}^{2}b_{k}+\kappa^{'} b_{0}\neq 0$
 with
$$\kappa=\frac{\disp\int_{\mathbb{H}^{1}}|t|^{\frac{\beta}{2}}\frac{1-||z|^{2}-it|^{2}}
{\Big|1+|z|^{2}-it\Big|^{6}}\theta_{0}\wedge
d\theta_{0}}{\disp\int_{\mathbb{H}^{1}}|x_{1}|^{\beta}
\frac{1-||z|^{2}-it|^{2}}
{\Big|1+|z|^{2}-it\Big|^{6}}\theta_{0}\wedge
d\theta_{0}},\quad\quad\quad\quad\kappa^{'}=\frac{\disp\int_{\mathbb{H}^{1}}\frac{|t|^{\frac{\beta}{2}}}
{\Big|1+|z|^{2}-it\Big|^{4}}\theta_{0}\wedge
d\theta_{0}}{\disp\int_{\mathbb{H}^{1}} \frac{|x_{1}|^{\beta}}
{\Big|1+|z|^{2}-it\Big|^{4}}\theta_{0}\wedge d\theta_{0}}$$\\

The function $\overset{[\beta]}{\underset{p=0}{\sum}}\big|\nabla^{p}
\mathcal{R}( x ) \big|\; \|x\|_{\mathbb{H}^{1}}^{-\beta-r} = o(1)$
as $x$ approaches $\xi_{i}$, $\nabla^{r} $ denotes  all possible
partial
derivatives of order $r$ and $[\beta]$ the integer part of $\beta.$\\\mbox{}\\
In this work, we will focus on the case where a collection of the
critical points of $K$ satisfy $ \beta=\beta(\xi_i)=2.$ This case
was not covered in the results of \cite{Gamara Riahi impact,Gamara
Riahi interplay,Gamara Riahi multiplicity}. So, here we suppose  $2
\leq \beta=\beta(\xi_i) < 4.$ Let
$$\begin{array}{l}
\mathcal K_{1} := \Big\{\xi_{i} \in\mathcal K\ \text{ such that}
\eqref{3'}\ \text{  is satisfied with}\; \beta=\beta(\xi_i)=2\ \
\text{and}\
 \sum_{k=1}^{2} b_{k}+\kappa^{'} b_{0}< 0\Big\}\\\\
\mathcal K_{2} := \Big\{ \xi_{i} \in\mathcal K\ \text{ such that}\
\eqref{3'}\ \text{  is satisfied with}\; \beta=\beta(\xi_i)>2 \ \
\text{and}\; \sum_{k=1}^{2} b_{k}+\kappa^{'} b_{0}< 0 \Big\}.
\end{array}$$

The index of the function $K$ at $\xi_{i}\in \mathcal K,$ denoted by
$m(\xi_i)$, is the number of strictly negative coefficients
$b_k(\xi_i)$:
$$m(\xi_i)=\#\Big\{b_{k}(\xi_i);b_{k}(\xi_i)<0\Big\}.$$
 For each  p-tuple $(\xi_{i_{1}},...,\xi_{i_{p}})\in (\mathcal K_1)^{p}$ ( $\xi_{i_{l}}\neq \xi_{i_{j}}$
if $l\neq j$), we associate the  matrix
$M(\xi_{i_{1}},...,\xi_{i_{p}})=(M_{st})_{1\leq s,t\leq p}$
\begin{eqnarray}\label{4}
 \begin{array}{rcl}
M_{ss} &=& \displaystyle{- c \frac{\sum_{k=1}^{2} b_{k}+\kappa' b_{0}}{2K^{2}(\xi_s)}}\\
M_{st} &=& \displaystyle{- c'
\frac{G(\xi_s,\xi_t)}{[K(\xi_s)K(\xi_t)]^{\frac{1}{2}}},\quad {\rm
for}\; s\neq t}\end{array}
\end{eqnarray}
where $c=
\int_{\mathbb{H}^{1}}\frac{|x_{1}|^{2}}{|1+|z|^{2}-it|^{4}}$ and
$c'= 2\pi \omega_{3},$ $\omega_{3}$ is the volume of the unit
Koranyi's ball.

 We say that $K$ satisfies  condition
$\mathbf{(C)}$ if:
\begin{eqnarray}\label{5'}
 \begin{array}{r}\hbox{{\rm  for each p-tuple }}  (\xi_{i_{1}},...,\xi_{i_{p}})\in (\mathcal K_1)^p\; \hbox{{\rm the corresponding matrix}}\; (M_{st})
 \;\hbox{{\rm is non degenerate}}.\end{array}
 \end{eqnarray}
In this case, we denote by $\varrho(\xi_{i_{1}},...,\xi_{i_{p}})$ the least eigenvalue of the matrix  $M(\xi_{i_{1}},...,\xi_{i_{p}}).$\\

 Next, we define the sets
 \begin{center}
 $\bullet$ $\mathcal K_{1}^{+}:=\bigcup_{p}\left\{  (\xi_{i_{1}},...,\xi_{i_{p}})\in( \mathcal K_{1})^{p},\
  \varrho(\xi_{i_{1}},...,\xi_{i_{p}})>0\right\}$
\end{center}
and
 \vspace*{5.mm}
 \begin{center} $\bullet$ $l^{+}:=\max\left\{p\in\mathbb{N} \ \textrm{ s.t }\exists\ (\xi_{i_{1}},...,\xi_{i_{p}})\in \mathcal K_{1}^{+}\right\}.$\vspace{,3cm}
\end{center}
For $(\xi_{i_{1}},...,\xi_{i_{p}})\in \mathcal K_{1}^{+},$ let
$i(\xi_{i_{1}},...,\xi_{i_{p}}):=4p-1- \sum_{j=1}^p
m(\xi_{i_{j}})$\\
and
\begin{eqnarray}\label{L0}
L_{0}:=\max\left\{ \Big\{i(\xi_{i_{1}},...,\xi_{i_{p}});\
(\xi_{i_{1}},...,\xi_{i_{p}})\in \mathcal K_{1}^{+}\}\cup
\Big\{3-m(\xi);\ \xi\in \mathcal K_{2}\Big\}\right\}
\end{eqnarray}
%%%%%%%%%%%%%%%%%%%%%%%%%%%%%%%%%%%%%%%%%%%%%%%%%%%%%%%%%%%%%%%%%%%%%%%%%%%%%
\par The main results of this paper are
\begin{thm} \label{existence}
Let $K$ be  a $C^{2}$ positive function on $\mathbb{S}^{3}$
satisfying  the \textbf{$\beta-$flatness condition} and condition
$\mathbf{(C)},$  if there exists a positive integer $k$ such that:

\begin{enumerate}
\item
$$\displaystyle\sum_{\scriptsize\begin{array}{l}\xi\in \mathcal K_{2}\\ m(\xi)\geq 4-k\end{array}}(-1)^{m(\xi)+1}+\;\;
\sum_{p=1}^{l^{+}}\sum_{\scriptsize{\begin{array}{l}(\xi_{i_{1}},.,\xi_{i_{p}})\in \mathcal K_{1}^{+}\\
 \sum_{j=1}^p
m(\xi_{i_{j}})\geq 4p-k \end{array}}}
(-1)^{i(\xi_{i_{1}},...,\xi_{i_{p}})}\neq 1$$

%\vspace{,3cm}
\item
$\forall (\xi_{i_{1}},...,\xi_{i_{p}})\in \mathcal K_{1}^{+}, \;
\sum_{j=1}^{p} m(\xi_{i_{j}})\neq 4p -(k+1)$\; and\; $\forall
\xi_i\in\mathcal K_2,$\; $3-m(\xi_i)\neq k.$

\end{enumerate}

Then, there exists a solution $\omega$ to the problem $\eqref{1}$
such that
$$m(\omega)\leq k,$$
where $m(\omega)$ is the Morse index of $\omega$, defined as the
dimension of the space of negativity
of the linearized operator $\mathcal{L}(\delta ) := L_{\theta}(\delta ) -3 \omega^{2}\delta. $\\
\end{thm}
%%%%%%%%%%%%%%%%%%%%%%%%%%%%%%%%%%%%%%%%%%%%%%%%%%%%%%%%%%%%%%%%%%%%%%%%%%%%%%%%%%%%%%%%%%%%%%
Under the hypothesis of Theorem \ref{existence}, if we denote
$\mathcal{S}_{k}$  the set of solutions of $\eqref{1}$ having their
Morse indices  less than or equal to $k$. We have
\begin{thm}\label{multiplicity}
$$\#\mathcal{S}_{k}\geq\left|1+\sum_{\scriptsize\begin{array}{l}\xi\in \mathcal K_{2}\\ m(\xi)\leq 4-k \end{array}}(-1)^{m(\xi)}-\;
\sum_{p=1}^{l^{+}}\sum_{\scriptsize\begin{array}{l}(\xi_{i_{1}},.,\xi_{i_{p}})\in \mathcal K_{1}^{+}\\
\sum_{j=1}^p m(\xi_{i_{j}})\geq 4p-k \end{array}} (-1)^{\sum_{j=1}^p
m(\xi_{i_{j}})}\right|$$

\end{thm}

%%%%%%%%%%%%%%%%%%%%%%%%%%%%%%%%%%%%%%%%%%%%%%%%%%%%%%%%%%%%%%%%%%55
The proofs of Theorems $\ref{existence}$ and $\ref{multiplicity}$
will be obtained by a contradiction argument. Therefore, we assume
that equation $\eqref{1}$ has no solution.
 Our approach  involves a Morse lemma at infinity, it relies on the construction of a suitable pseudo gradient for the functional $J.$
  The Palais-Smale condition is satisfied along the decreasing flow lines of this pseudo gradient, as long as these
 flow lines do not enter the neighborhood of a finite number of critical points of $K$ where the related matrix given in $(\ref{4})$
  is positive definite.\\
 This paper is organized as follows: in section 2, we recall the local structure
  of the Heisenberg group, the extremals for the Yamabe functional on $\mathbb H^1$ and the Cayley transform. In section 3,
  we give the expansion of the
new functional $J$ near its  critical points at infinity.
 Section 4 is devoted to the construction of
a Morse Lemma at infinity for the functional $J$. The Morse lemma
 is based on the construction  of a pseudo gradient for
$J$  near its  critical points at infinity, using an appropriate
change of variables. The proofs of our main results, Theorems
\ref{existence}
 and $\ref{multiplicity}$ will be the purpose of section 5.
The last section is an appendix, where some technical estimates
 are given.

\section{Preliminary Tools:}
\hspace*{0.5cm}The Heinserberg group $\mathbb{H}^{1}$ is the Lie
group whose underlying manifold is ${\mathbb{C}}\times
{\mathbb{R}}$, with coordinates $g=(z,t)$ and
 group law  given by: $g\cdot g^{\prime } = (z,t)\cdot
(z^{\prime },t^{\prime }) = (z+z^{\prime }, t + t^{\prime} + 2
Im\,z.\bar{z}^{\prime })$. We define a norm in $\mathbb{H}^{1}$ by
$\| g \|_{\mathbb{H}^{1}} = \| (z,t) \|_{\mathbb{H}^{1}} = (\| z
\|^{4} + t^{2})^{\frac{1}{4}}$, and dilations by $g =
(z,t)\rightarrow \lambda g = (\lambda z, \lambda^{2}t)$, $\lambda >
0$. The Cauchy Riemann structure on $\mathbb{H}^{1}$ is given by the
left invariant vectors fields: $Z=\frac{\partial }{\partial z} +
i\bar{z}\frac{\partial }{\partial t}$, $\bar{Z} = \frac{\partial
}{\partial \bar{z}} - iz\frac{\partial }{\partial t}$, which are
homogenous of degree $-1$ with respect to the dilations, the
associated  contact form is $\theta_{0} = dt + i(z d\bar{z} -
\bar{z} dz)$. We denote by $\Delta_{\theta _{0}}$ the sublaplacian
operator, $\Delta_{\theta_{0}} = -\;\frac{1}{2}(Z\bar{Z} +
\bar{Z}Z)$ and since the Webster scalar curvature $R_{\theta _{0}}$
is zero, the conformal laplacian $L_{0}$ is a multiple of the
sublaplacian operator,
$L_{0}=(2+\frac{2}{n})\Delta_{\theta _{0}}$.\\
In \cite{J-L1}, Jerison and Lee showed that all solutions of
$\eqref{2}$ are obtained from
\begin{eqnarray*}
w_{(0,1)}(z,t)=\frac{c_{0}}{|1+|z|^{2}-it|}\;,\;\;c_{0} > 0,
\end{eqnarray*}
by left translations and dilatations on $\mathbb{H}^{1}$. That is
for $g_{0} = (z_{0},t_{0})$, $g = (z,t)$ in $\mathbb{H}^{1}$ and
$\lambda > 0$, we have
\begin{eqnarray*}
w _{(g_{0},\lambda )}(z,t)= c_{0}\frac{\lambda}{|1+\lambda
^{2}|z-z_{0}|^{2}-i\lambda ^{2}(t-t_{0}-2Im\, z_{0}\overline{z})|}
\end{eqnarray*}
%%%%%%%%%%%%%%%%%%%%%%%%%%%%%%%%%%%%%%%%%%%%%%%%%%%%%%%%%%%%%%%%%%%%%%%%%%%%%%%%%%%%%%%
Next, we will introduce the Cayley transform. Let $B^{2 }= \big\{z
\in \mathbb{C}^{2}\;/\; |z| < 1 \big\}$ be the unit ball in
$\mathbb{C}^{2}$ and $\mathcal{D}_{2} = \big\{ (z,w) \in
\mathbb{C}\times \mathbb{C}\;/\; Im(w) > | z |^{2}\big\}$ be the
Siegel domain. The boundary of the Siegel domain is:
$\partial\mathcal{D}_{2} = \big\{ (z,w) \in \mathbb{C}\times \mathbb{C}\;/\; Im(w) = | z |^{2}\big\}$.\\
\begin{defn}\cite{D-T} The Cayley transform is the
correspondence between the  unit ball $B^{2}$ in $\mathbb{C}^{2 }$
and the Siegel domain $\mathcal{D}_{2},$ given by
\begin{eqnarray*}
\mathcal{C}(\zeta ) = \Big( \frac{\zeta_{1}}{ 1 +
\zeta_{2}}\,\,,\,\,\, i\,\,\frac{1 - \zeta_{2}}{ 1 +
\zeta_{2}}\Big)\,\,;\,\quad \zeta = (\zeta_{1}, \zeta_{2})\,\,,\quad
1 + \zeta_{2} \neq 0
\end{eqnarray*}
\end{defn}
The Cayley transform gives a biholomorphism of the unit ball $B^{2}$
in $\mathbb{C}^{2 }$ onto the Siegel domain $\mathcal{D}_{2}$.
Moreover, when restricted to the sphere minus a point, $\mathcal{C}$
gives a $CR$ diffeomorphism.
\begin{eqnarray*}
\mathcal{C }: \mathbb{S}^{3}\backslash {(0,-1)}\longrightarrow
\partial\mathcal{D}_{2}.
\end{eqnarray*}
Let us recall the CR diffeomorphism
\begin{eqnarray*}
\begin{array}{cccl}
 f \,: & \mathbb{H}^{1} & \longrightarrow  & \partial\mathcal{D}_{2} \\
   & (z, t) & \longmapsto & f(z,t) = (z, t + i |z|^{2})\,
\end{array}
\end{eqnarray*}
with the obvious inverse $f^{-1}(z,w) = (z, Re(w))$, $z\in
\mathbb{C}$, $ w \in \mathbb{C}$. We obtain the $CR$ equivalence via
this mapping:
\begin{eqnarray*}
\begin{array}{cccl}
F\,: & \mathbb{S}^{3}\backslash {(0,-1)} & \longrightarrow  & \mathbb{H}^{1} \\
   & \zeta = (\zeta_{1}, \zeta_{2}) & \longmapsto & (z,t)
   = \big(\frac{\zeta_{1}}{ 1 + \zeta_{2}},
   i\,\,\frac{2 Im\zeta_{2}}{| 1 + \zeta_{2}|^{2}}\big)
\end{array}
\end{eqnarray*}
with inverse
\begin{eqnarray*}
\begin{array}{cccl}
 F^{-1}\,: &  \mathbb{H}^{1} & \longrightarrow & \mathbb{S}^{3}\backslash {(0,-1)} \\
    & (z,t)  & \longmapsto & \zeta = \big(\frac{2 z}{ 1 + |z|^{2} - it},\,i\,\frac{1 - |z|^{2} + it }{
1 + |z|^{2} - it }\big).
\end{array}
\end{eqnarray*}
With the following choice of  contact form on $\mathbb{S}^{3}$ ( the
standard one)
\begin{eqnarray*}
\theta_{1} = i \overset{2}{\underset{j=1}{\sum}}\big( \zeta_{j}d
\overline{\zeta}_{j} - \overline{\zeta}_{j} d \zeta_{j} \big).
\end{eqnarray*}
We obtain
$F^{*}(4 (c_{0}^{-1}w_{(0,1)})^{2}\theta_{0}) = \theta_{1}$.\\
Let us differentiate and take into account that $w_{(0,1)}(F(\zeta))
= c_{0}|1 + \zeta_{2}|^{2}$, we obtain
\begin{eqnarray*}
d\theta_{1} =\big(\frac{d\zeta_{2}}{ 1 + \zeta_{2}} +
\frac{d\bar{\zeta}_{2}}{ 1 + \bar{\zeta}_{2} }  \big)\wedge
\theta_{1} + |1 + \zeta_{2}|^{2}F^{*}( d \theta_{0})
\end{eqnarray*}
and
\begin{eqnarray*}
\theta_{1}  \wedge d \theta_{1} = |1 +
\zeta_{2}|^{4}F^{*}(\theta_{0} \wedge d \theta_{0})
\end{eqnarray*}
We introduce the following function for each $(\zeta_{0}, \lambda)$
on $\mathbb{S}^{3}\times]0, + \infty[$
\begin{eqnarray}\label{dela0}
\delta_{(\zeta_{0},\lambda )}(\zeta )=|1+\zeta_{2}|^{-1}
w_{(F(\zeta_{0}),\lambda )} \circ F(\zeta )
\end{eqnarray}
 We have
$L_{\theta_{1}}\delta_{(\zeta_{0},\lambda
)}=\delta_{(\zeta_{0},\lambda )}^{3},$ i.e
 $\delta_{(\zeta_{0},\lambda )}$ is a solution of the Yamabe problem on $\mathbb{S}^{3}$.\\
We also have
\begin{eqnarray}\label{S-H-1}
\int_{\mathbb{S}^{3}}L_{\theta_{1}}\delta_{(\zeta_{0},\lambda )}
\,\,\delta_{(\zeta_{0},\lambda )}
 \,\,\theta_{1}\wedge d\theta_{1} =
\int_{\mathbb{H}^{1}}L_{\theta_{0}}w_{(g_{0},\lambda )}\,\,
w_{(g_{0},\lambda )}\,\,\theta_{0}\wedge d\theta_{0}\;
\end{eqnarray}
and
\begin{eqnarray}\label{S-H-2}
\int_{\mathbb{S}^{3}}|\delta_{(\zeta_{0},\lambda )}|^{ 4}\theta_{1}
\wedge d\theta_{1} = \int_{\mathbb{H}^{1}}| w_{(g_{0},\lambda
)}|^{4} \;\theta_{0}\wedge d\theta_{0},
\end{eqnarray}
where \; $g_{0} = F(\zeta_{0})$,\;and  $g = F(\zeta)$.\\
 As a consequence, the variational
formulation for $\eqref{1}$ is equivalent
to the one for $\eqref{2}$.\\
%%%%%%%%%%%%%%%%%%%%%%%%%%%%%%%%%%%%%%%%%%%%%%%%%%%%%%%%%%%%%%%%%%%%%%%%%%%%%%%%%%%%%%%%%%%%%%%%%%%%%%%%%%%
%%%%%%%%%%%%%%%%%%%%%%%%%%%%%%%%%%%%%%%%%%%%%%%%%%%%%%%%%%%%%%%%%%%%%%%%%%%%%%%%%%%%%%%%%%%%%%%%%%%%%%%%%%%
\subsection{ Cauchy Riemann Functional}

Problem $\eqref{1}$ has a nice variational structure, with
associated Euler functional:
$$J(u)=\disp\frac{\int_{\mathbb{S}^{3}} L_{\theta_{1}}u\; u\; \theta_{1}  \wedge d\theta_{1}}
{(\int_{\mathbb{S}^{3}} K\;u^{4}\; \theta_{1} \wedge
d\theta_{1})^{\frac{1}{2}}}, \quad u\in S_1^2(\mathbb{S}^{3})$$
where $S_1^2(\mathbb{S}^{3})$ is the completion of
$C^{\infty}(\mathbb{S}^{3})$ by means of the norm
$\left\|u\right\|^{2}=\int_{\mathbb{S}^{3}}
L_{\theta_{1}}u\,u\,\theta_{1}\wedge d\theta_{1}.$
\\ Let $\quad \Sigma =
\big\{u\in S_{1}^{2}(\mathbb{S}^{3}) / \left\| u\right\| = 1 \big\}$
and
$\quad\Sigma^{+}=\big\{u\in \Sigma / \;u\geq 0 \big\}$. \\

 The functional $J$ fails to satisfy the Palais-Smale condition denoted by $(P.S)$ on $\Sigma^+$,
that is: there exist noncompact sequences along which the functional
$J$ is bounded and its gradient goes to zero.  %The failure of  the
%(P.S) condition% has been analyzed for the Riemannian case throughout the works of
% \cite{Bahri Critical,Bahri invariant,B C H,B C C H,K W,Li,S U,St}. For the Cauchy Riemannian  case,
  A complete description
  of sequences failing to satisfy (P.S) is given  in \cite{Gamara Yacoub yamabe conjecture}.
 A solution  $u$ of $\eqref{1}$ is a critical point of $J$ subject to the constraint
$u\in \Sigma^+.$

\subsection{Characterization of the sequences failing to satisfy the (P.S) condition}

\par In the case we study, we have  the presence of multiple blow-up points.
 We begin by defining  the sets of potential critical points at infinity of
the functional $J.$\\
For any $\varepsilon >0$ and $p\in \mathbb{N}^{+}$, let:

$$\begin{array}{rlc}\label{critinf}
V(p,\varepsilon )&=&\left\lbrace\begin{array}{l}u\in
\Sigma^{+};\;\exists \;(a_{1},\ldots,a_{p})\in  \mathbb{S}^{3},\;
\alpha_{1},\ldots,\alpha_{p}>0
                    \;{\rm and}\;  (\lambda_{1},\ldots,\lambda_{p})\in (\varepsilon ^{-1},\infty )^{p}\;{\rm  such \;that}\\\\
                    \left\|u-\displaystyle{\sum_{i=1}^{p}\frac{\alpha_i
                    \delta _{a_{i}, \lambda_{i}}}{K(a_{i})^{\frac{1}{2}}}}\right\|_{S_{1}^{2}(\mathbb{S}^{3})} <
                    \varepsilon,\;
                    \varepsilon_{ij} < \varepsilon,\;| \displaystyle{\frac{\alpha_{i}^{2}K(a_{i})}{\alpha_{j}^{2}K(a_{j})}} - 1| < \varepsilon,
                     \;\;\forall\;1\leq i\neq j\leq p\\
                      \varepsilon_{ij} = \big( \displaystyle{\frac{\lambda_{i}}{\lambda_{j}}+ \frac{\lambda_{j}}{\lambda_{i} }} + \lambda_{i}\lambda_{j}
 (d(a_{i}, a_{j})^{2} \big)^{-1}. \end{array}\right\rbrace
\end{array}$$

For $\omega$ a solution of $\eqref{1}$ we also define the set
\begin{eqnarray}\label{critinfinisol}
V(p,\varepsilon,\omega) =\Big\{u\in\Sigma^+;\exists\ \alpha_{0}>0\;
/\;\ u-\alpha_{0}\omega\in V(p,\varepsilon)\ \rm{and}\
|\alpha_{0}^{2}J(u)^{2}-1|<\varepsilon\Big\}.
\end{eqnarray}

 We then proceed as in
\cite{Gamara Yacoub yamabe conjecture} Proposition $8$ to
characterize the sequences which violate the (P.S) condition as
follows:
\begin{prop}(\cite{Gamara Yacoub yamabe conjecture}) Let $\{u_{k}\}$ be a sequence such that $\partial J(u_{k})\rightarrow 0$
and $J(u_{k})$ is bounded.  There exist an integer $p\in
\mathbb{N}^{*}$, a sequence $\varepsilon_{k}\rightarrow
0\;\;(\varepsilon_{k}>0)$ and an extracted subsequence of
$\{u_{k}\}$, again denoted by $\{u_{k}\}$, such that $u_{k}\in
V(p,\varepsilon_{k}).$
\end{prop}
Then, we consider the following minimization problem for a function
$u\in V(p,\varepsilon ),$ with $\varepsilon $ small
\begin{eqnarray}\label{minor1}
\min_{\alpha_i>0,\lambda_i>0,a_i\in
\mathbb{S}^{3}}\|u-\sum_{i=1}^{p}\alpha_{i}
\delta_{a_{i},\lambda_{i}}\|_{S_{1}^{2}(\mathbb{S}^{3})}
\end{eqnarray}
We obtain as showed in \cite{Bahri invariant} and \cite{Gamara
prescribed},  the following parametrization of the set $V(p,
\varepsilon):$
\begin{prop} (\cite{Gamara
prescribed}) For any $p\in {\Bbb N}^{*}$, there exists
$\varepsilon_{p}>0$ such that, for any $0 < \varepsilon  <
\varepsilon_{p}$, $u\in V(p,\varepsilon )$, the minimization problem
$\eqref{minor1}$ has a unique solution
$(\bar{\alpha}_{1},\ldots,\bar{\alpha}_{p},\bar{\lambda}_{1},\ldots,\bar{\lambda}_{p},\bar{a}_{1},\ldots,\bar{a}_{p})$
\big(up to permutation on the set of indices $\{1,\ldots,p\}$\big).
In particular, we can write $u \in V(p, \varepsilon)$ as follows
\end{prop}
 $u =  \overset{p}{\underset{i=1}{\sum}}\overline{\alpha}_{i}\delta _{\overline{a}_{i},\overline{\lambda}_{i}} + v$,
  where $v \in S_1^2(\mathbb{S}^{3})$ satisfies:
\begin{eqnarray*}
(V_{0})\,\,\left\{
\begin{array}{ll}
\langle v,\delta_{a_{i},\lambda_{i}}\rangle _{S_1^2(\mathbb{S}^{3})} & =0 \\
\langle v,\frac{\partial \delta_{a_{i},\lambda_{i}}}{\partial
a_{i}}\rangle _{S_1^2(\mathbb{S}^{3})} &
= 0\qquad\quad i=1,2,...,p. \\
\langle v,\frac{\partial \delta_{a_{i},\lambda_{i}}}{\partial
\lambda_{i}}\rangle _{S_1^2(\mathbb{S}^{3})} & = 0.
\end{array}
\right.
\end{eqnarray*}
Here $< , >$ denotes the $L$-scalar product defined on
$S_{1}^{2}(\mathbb{S}^{3})$ by
\begin{eqnarray}\label{8}
<u,v>=\int_{\mathbb{S}^{3}}L_{\theta_{1}}u v  \;\theta_{1}\wedge
d\theta_{1}.
\end{eqnarray}

Next, we will focus on  the behavior of the functional $J$ with
respect to the variable $v$. We will prove the existence of a unique
$\bar{v}$ which minimizes
$J(\overset{p}{\underset{i=1}{\sum}}\alpha_{i}\delta
_{a_{i},\lambda_{i}}+v)$ with respect to $v\in H_{\varepsilon
}^{p}(a,\lambda )$, where
\begin{eqnarray*}
H_{\varepsilon }^{p}(a,\lambda ) = H_{\varepsilon }^{p} (\delta
_{a_{1},\lambda_{1}},\ldots,\delta _{a_{p},\lambda_{p}}) =
\Big\{v\in S_{1}^{2}(M)\text{ \ }/v\text{ satisfies }(V_{0}) \text{
and }\|v\|<\frac{\varepsilon }{p}\Big\}.
\end{eqnarray*}
\begin{prop}\cite{Gamara prescribed} There exists a $C^{1}$-map which associates to each $u\in V(p,\varepsilon ),$
$\varepsilon$ small,  $\bar{v}=\bar{v}(\alpha, a, \lambda )$ such
that $\bar{v}$ is unique and minimizes
$J(\overset{p}{\underset{i=1}{\sum}}\alpha_{i}\delta _{a_{i},
\lambda_{i}} + v)$, with respect to $v\in H_{\varepsilon
}^{p}(a,\lambda ).$   We have the following estimate

$$\|\bar{v }\|  \leq c_1\left(\displaystyle{
\underset{i\leq p}{\sum} (\frac{|\nabla K(a_i)|}{\lambda_i}+
\frac{1}{\lambda_i^2})+\underset{k\neq
r}{\sum}\varepsilon_{kr}\sqrt{{\rm Log}(\varepsilon_{kr}^{-1})}}
\right)$$$\Box$
\end{prop}
For  $\omega$  a solution of $\eqref{1}$, we obtain a
parametrization of the set $V(p,\varepsilon,\omega)$ as follows

\begin{prop}\label{minimum} There is $\varepsilon_0>0$ such that if $\varepsilon\leq\varepsilon_0$ and $u\in V(p,\varepsilon,\omega),$
  the problem
$$\underset{\alpha_i>0,\;\lambda_i>0,\; a_i \in \mathbb{S}^{3},\; h\in T_{\omega}(W_u(\omega))}{\min}\|u-\sum_{i=1}^p\alpha_{i}\delta _{a_{i},\lambda_i}- \alpha_0(\omega + h)\|$$
has a unique solution $(\overline \alpha, \overline
\lambda,\;\overline a,\;\overline h).$ Thus, we write $u$ as:
$$u =  \sum_{i=1}^p\overline\alpha_{i}\delta _{(\overline a_{i},\overline \lambda_i)}+ \overline\alpha_0(\omega +\overline h) + v,$$
where $v$ belongs to $S_1^2(\mathbb{S}^{3})\cap
T_{\omega}(W_s(\omega))$ and satisfies
$(V_0),\;T_{\omega}(W_u(\omega))$ and $T_{\omega}(W_s(\omega))$ are
respectively, the tangent spaces at $\omega$ to the unstable and
stable manifolds of $\omega.$
\end{prop}
{\bf Proof:} The proof is similar to the one given in \cite{Bahri
invariant}.

\section{Asymptotic Analysis of the Functional}
\subsection{Domination Property:  Hierarchy of the Critical point at
infinity } We first introduce some definitions and notations due to
Bahri \cite{ Bahri Critical, Bahri invariant}. Let $\partial J$
denotes the gradient of the functional $J.$
\begin{defn}
A critical point at infinity of $J$ on $\Sigma^+$ is a limit of a
flow line $u(s)$ of the equation:
$$\left\lbrace
\begin{array}{ll}
\displaystyle{\frac{\partial u }{\partial s}} & =-\partial J(u)\vspace*{,1cm} \\\
u(0) & =u_{0}
\end{array}
\right.$$ such that $u(s)$ remains in $V(p,\varepsilon(s),\omega)$
for $s\geq s_0,$ $\omega$ is zero or a solution of (\ref{existence})
and $\varepsilon(s)$ satisfies
$\underset{s\longrightarrow\infty}{\lim}\varepsilon(s)= 0.$
  One can write $u(s) =
\sum_{i=1}^p\alpha_{i}(s)\delta _{(a_{i}(s), \lambda_i(s))}+
\alpha_{0}(s)(\omega + h(s)) + v(s).$ Let
$a_i:=\underset{s\longrightarrow\infty}{\lim}a_i(s)$ and $\alpha_i:=
\underset{s\longrightarrow\infty}{\lim}\alpha_i(s),$ we denote such
a critical point at infinity by $$\xi_{\infty}\;{\rm
or}\;(a_1,\ldots,a_p)_{\infty}\;{\rm or}\;
\sum_{i=1}^p\alpha_i\delta _{(a_i,\infty)}\;{\rm or}\;
\sum_{i=1}^p\alpha_{i}(s)\delta _{(a_{i}, \infty)}+ \alpha_{0}
\omega.$$
\end{defn}
A critical point at infinity is called of $\omega-$type if
$\omega\neq 0.$\\
 As for a usual critical point, to a critical point at infinity $\xi_{\infty}$ are
associated stable and unstable manifolds which we denote by
$W_{s}(\xi_{\infty})$ and
 $W_{u}(\xi_{\infty}).$ These manifolds allow to compare critical points at infinity by what we call a "domination property",
  one can see \cite{Bahri invariant, Gamara prescribed}, where a detailed description of theses manifolds is given.

 \begin{defn}
A critical point at infinity $\xi_{\infty}$ is said to be dominated
by another critical point at infinity $\xi'_{\infty},$
  if $$W_{s}(\xi_{\infty}) \cap W_{u}(\xi_{\infty}^{'})\neq \emptyset$$ and we write  $\xi_{\infty^{'}} > \xi_{\infty}.$ $\Box$\\
\end{defn}
If we assume that the intersection $W_{s}(\xi_{\infty}) \cap
W_{u}(\xi_{\infty}^{'})$  is transverse, then we obtain

$\textrm{index} (\xi_{\infty}^{'})\geq \textrm{index} (\xi_{\infty})
+ 1.$

\subsection{Asymptotic Analysis of the functional in the set $V(p,\varepsilon,\omega),$ $\omega\neq 0$}
In this section, we expand the functional $J$ in
$V(p,\varepsilon,\omega),$ for $\omega$ a non null solution of
$\eqref{1}$ in the aim to detect  the critical points or critical points at infinity of $J$ in this set and we prove that:\\
 for any $p \in \mathbb N^{\ast}$, there are no critical point or critical point at infinity of $J$ in the set $V (p, \varepsilon,\omega)$.
More precisely, using Proposition \ref{minimum}, we will write $u
\in V (p, \varepsilon,\omega)$ as
$u=\disp\sum_{i=1}^{p}\alpha_{i}\delta
_{a_{i},\lambda_{i}}+\alpha_{0}(\omega+h)+v,$ one obtain the
following expansion of $J:$

\begin{prop}\label{expansionJ(w)}
There exists $\varepsilon_0 >0$ such that for any
$u=\sum_{i=1}^p\alpha_{i} \delta _{(a_{i}, \lambda_i )}+ \alpha_0
(\omega + h ) + v \in V(p,\varepsilon,\omega),
\varepsilon<\varepsilon_0$
\begin{eqnarray*}
J(u)&=&\frac{S\sum_{i=1}^{p}\alpha_{i}^{2}+\alpha_{0}^{2}\|w\|^{2}}
{(S\sum_{i=1}^{p}\alpha_{i}^{4}K(a_{i})+\alpha_{0}^{4}\|w\|^{2})^{\frac{1}{2}}}
\Big[1-\frac{c_{2}\alpha_{0}}{\gamma_{1}}\sum_{i=1}^{p}\alpha_{i}\frac{w(a_{i})}{\lambda_{i}}
-\frac{1}{\gamma_{1}}\sum_{i\neq j}\alpha_{i}\alpha_{j}c_{ij}\varepsilon_{ij}+f_{1}(v)\\
&&+Q_{1}(v,v) +f_{2}(h)+\alpha_{0}^{2}Q_{2}(h,h) +o\left(\sum_{i\neq
j}\varepsilon_{ij}+\sum_{i=1}^{p}\frac{1}{\lambda_{i}}+\|v\|^{2}+\|h\|^{2}\right)\Big]\end{eqnarray*}
\begin{eqnarray*}&&f_{1}(v)=-\frac{1}{\eta_{1}}
\int_{\mathbb{S}^{3}}K(\sum_{i=1}^{p}\alpha_{i}\delta_{a_{i},\lambda_{i}})^{3}v,\\
&&f_{2}(h)=\frac{\alpha_{0}}{\gamma_{1}}
\sum_{i=1}^{p}\alpha_{i}<\delta_{a_{i},\lambda_{i}},h>_{L_{\theta}}-
\frac{\alpha_{0}}{\eta_{1}}\int_{\mathbb{S}^{3}}(\sum_{i=1}^{p}\alpha_{i}\delta_{a_{i},\lambda_{i}}+\alpha_{0}w)^{3}h,\\
\end{eqnarray*}

\begin{eqnarray*}
&&Q_{1}(v,v)=\frac{\|v\|^{2}}{\gamma_{1}}- \frac{3}{\eta_{1}}
\int_{\mathbb{S}^{3}}K\left(\sum_{i=1}^{p}(\alpha_{i}\delta_{a_{i},\lambda_{i}})^{2}+
(\alpha_{0}w)^{2}\right)v^{2},\\
&&Q_{2}(h,h)=\frac{\|h\|^{2}}{\gamma_{1}}- \frac{3}{\eta_{1}}
\int_{\mathbb{S}^{3}}K(\alpha_{0}w)^{2}h^{2},\\
&&c_{2}=c_{0}^{3}\int_{\mathbb{H}^{1}}\frac{1}{\Big|1+|z|^{2}-it \Big|^{3}}\theta_{0}\wedge d\theta_{0},\quad  S=c_{0}^{4}\int_{\mathbb{H}^{1}}\frac{1}{\Big|1+|z|^{2}-it \Big|^{4}}\theta_{0}\wedge d\theta_{0},\\
&&\eta_{1}=S\sum_{i=1}^{p}\alpha_{i}^{4}K(a_{i})+\alpha_{0}^{4}\|w\|^{2},\quad
\gamma_{1}=S\sum_{i=1}^{p}\alpha_{i}^{2}+\alpha_{0}^{2}\|w\|^{2}
\end{eqnarray*}
and $c_{ij}$ are  bounded positive constants.
\end{prop}

\textbf{Proof:}  we need to estimate
$$N=\|u\|^{2}\quad and\quad D^{2}=\int_{\mathbb{S}^{3}}Ku^{4}\theta\wedge d\theta.$$
Expanding $N$, we get
$$N:=\sum_{i=1}^{p}\alpha_{i}^{2}\|\delta_{i}\|^{2}+\alpha_{i}\alpha_{0} <\delta_{i},w+h>_{L_{\theta}}+\alpha_{0}(\|h\|^{2}+\|w\|^{2})+\|v\|^{2}+\sum_{i\neq j}
\alpha_{i}\alpha_{j}<\delta_{i},\delta_{j}>_{L_{\theta}}.$$ It
follows from [15] and elementary computations that
\begin{eqnarray*}
\|\delta_{i}\|^{2}&=&S,\\
<\delta_{i},\delta_{j}>_{L_{\theta}}&=&c_{ij}\varepsilon_{ij}(1+o(1)),\ \textrm{for}\ i\neq j\ and \\
<\delta_{i},w>_{L_{\theta}}&=&c_{2}\frac{w(a_{i})}{\lambda_{i}}+o(\frac{1}{\lambda_{i}}).
\end{eqnarray*}
Therefore
\begin{eqnarray*}
N=\gamma_{1}+2\alpha_{0}c_{2}\sum_{i=1}^{p}\alpha_{i}\frac{w(a_{i})}{\lambda_{i}}
+\alpha_{i}<\delta_{i},h>_{L_{\theta}} +\sum_{i\neq
j}\alpha_{i}\alpha_{j}+
\alpha_{0}^{2}\|h\|^{2}+\|v\|^{2}+o\left(\sum_{i=1}^{p}\frac{1}{\lambda_{i}}+
\sum_{i\neq j}\varepsilon_{ij}\right).
\end{eqnarray*}
For the denominator $D$, we compute it as follows
\begin{eqnarray*}
D^{2}&=&\int_{\mathbb{S}^{3}}K\left(\sum_{i=1}^{p}\alpha_{i}\delta_{i}\right)^{4}\theta\wedge d\theta+\int_{\mathbb{S}^{3}}K(\alpha_{0}w)^{4}\theta\wedge d\theta+4\alpha_{0}\int_{\mathbb{S}^{3}}K\left(\sum_{i=1}^{p}\alpha_{i}\delta_{i}\right)^{3}w\theta\wedge d\theta\\
&&+4\alpha_{0}^{3}\int_{\mathbb{S}^{3}}K\left(\sum_{i=1}^{p}\alpha_{i}\delta_{i}\right)w^{3}\theta\wedge
d\theta+4\int_{\mathbb{S}^{3}}K\left(\sum_{i=1}^{p}\alpha_{i}\delta_{i}+\alpha_{0}w\right)^{3}(\alpha_{0}h+v)
\theta\wedge d\theta\\
&&+12\int_{\mathbb{S}^{3}}K\left(\sum_{i=1}^{p}\alpha_{i}\delta_{i}+\alpha_{0}w\right)^{2}(\alpha_{0}^{2}h^{2}
+v^{2}+2hv)\theta\wedge d\theta
+O\left(\sum_{i=1}^{p}\int_{\mathbb{S}^{3}}w^{2}\alpha_{i}^{2}\delta_{i}^{2}\right)\\
&&+O(\|v\|^{3}+\|h\|^{3}).
\end{eqnarray*}
Where
\begin{eqnarray*}
&&\int_{\mathbb{S}^{3}}K\left(\sum_{i=1}^{p}\alpha_{i}\delta_{i}\right)^{4}\theta\wedge d\theta=\sum_{i=1}^{p}\alpha_{i}^{4}K(a_{i})S+4\sum_{i\neq j}\alpha_{i}^{3}\alpha_{j}K(a_{i})c_{ij}\varepsilon_{ij}+O\left(\sum_{i=1}^{p}\frac{1}{\lambda_{i}^{2}}\right),\\
&&\int_{\mathbb{S}^{3}}Kw^{4}\theta\wedge d\theta=\|w\|^{2},\quad \int_{\mathbb{S}^{3}}Kw^{3}\delta_{i}\theta\wedge d\theta=c_{2}\frac{w(a_{i})}{\lambda_{i}}+o(\frac{1}{\lambda_{i}}),\\
&&\int_{\mathbb{S}^{3}}(w^{2}\alpha_{i}^{2}\delta_{i}^{2}+w^{2}\alpha_{0}^{2}\delta_{i}^{2})\theta\wedge d\theta=o(\frac{1}{\lambda_{i}}),\\
&&\int_{\mathbb{S}^{3}}\left(\sum_{i=1}^{p}\alpha_{i}\delta_{i}+\alpha_{0}w\right)^{2}hv\theta\wedge
d\theta=
O\left(\int_{\mathbb{S}^{3}}\Big(\sum_{i=1}^{p}\delta_{i}^{2}+w^{-1}\sum_{i=1}^{p}\delta_{i}\Big)|h||v|\theta\wedge d\theta\right)\\
&&\hspace{5cm}=O\left(\|h\|^{3}+\|v\|^{3}+\sum_{i=1}^{p}\frac{1}{\lambda_{i}^{3}}\right)
\end{eqnarray*}
where we have used that $v\in T_{w}(W_{s}(w))$ and $h$ in
$T_{w}(W_{u}(w))$.\\
Next, we focus on the linear form in  $v\in T_{w}(W_{s}(w))$, we
obtain
\begin{eqnarray*}
\int_{\mathbb{S}^{3}}K\left(\sum_{i=1}^{p}\alpha_{i}\delta_{i}+\alpha_{0}w\right)^{3}v\theta\wedge
d\theta&=&
\int_{\mathbb{S}^{3}}K\left(\sum_{i=1}^{p}\alpha_{i}\delta_{i}\right)^{3}v\theta\wedge
d\theta
+O\left(\sum_{i=1}^{p}\int_{\mathbb{S}^{3}}(\alpha_{i}^{2}\alpha_{0}\delta_{i}^{2}w+
\alpha_{i}\alpha_{0}^{2}\delta_{i}w^{2})|v|\right)\\
&=&f_{1}(v)+O\left(\sum_{i=1}^{p}\frac{\|v\|}{\lambda_{i}}\right).
\end{eqnarray*}
Finally, for the partial quadratic forms in $v$ and $h,$ we obtain
\begin{eqnarray*}
&&\int_{\mathbb{S}^{3}}K\left(\sum_{i=1}^{p}\alpha_{i}\delta_{i}+\alpha_{0}w\right)^{2}h^{2}\theta\wedge d\theta=\alpha_{0}^{2}\int_{\mathbb{S}^{3}}Kw^{2}h^{2}\theta\wedge d\theta+o(\|h\|^{2})\\
&&\int_{\mathbb{S}^{3}}K\left(\sum_{i=1}^{p}\alpha_{i}\delta_{i}+\alpha_{0}w\right)^{2}v^{2}\theta\wedge
d\theta=\sum_{i=1}^{p}\int_{\mathbb{S}^{3}}K(\alpha_{i}\delta_{i})^{2}v^{2}\theta\wedge
d\theta+ \alpha_{0}^{2}\int_{\mathbb{S}^{3}}Kw^{2}v^{2}\theta\wedge
d\theta+o(\|v\|^{2})
\end{eqnarray*}
Combining these results and the fact that $\frac{\alpha_{i}^{2}K(a_{i})}{\alpha_{j}^{2}K(a_{j})}=1+o(1)$ the proposition follows. $\Box$\\
Next, we state the following result
 \begin{lem} \cite{Bahri invariant}
\begin{itemize}
\item[a-] $Q_1(v,v)$ is a positive definite  quadratic form   on
$$E_v =\left\{v\in S_1^2(M) \;{\hbox{\rm such that}}\; v\in T_{\omega}(W_s(\omega)) \;{\hbox{\rm and}}\;v \;{\hbox{\rm satisfies}}(V_0)\;\right\}.$$
\item[b-] $Q_2(h,h)$ is a negative definite quadratic form  on $T_{\omega}(W_u(\omega)).$
\end{itemize}
\end{lem}
{\bf Proof:} The proof of this lemma is similar to the one given in
\cite{Bahri invariant},  for more details one can see the appendix
of
\cite{Gamara prescribed}, where necessary modifications are given.$\Box$\\

Using the Lemma above one can perform the expansion of the
functional $J$ given in Proposition  \ref{expansionJ(w)} after an
adequate change of variables. More precisely,  we obtain
\begin{prop}
Let $u  =  \sum_{i=1}^p\alpha_{i} \hat\delta _{(a_{i} , \lambda_i)}+
\alpha_0(\omega + h) + v \in V(p,\varepsilon,\omega).$ There is an
optimal $(\overline v,\overline h)$ and a change of variables
$v-\overline v\rightarrow V$ and $h-\overline h\rightarrow H$ such
that
$$
J(u)= J\left(\sum_{i=1}^p\alpha_{i} \hat\delta _{(a_{i} ,
\lambda_i)}+ \alpha_0\omega + \overline h+\overline v\right) +
\|V\|^2- \|H\|^2.
$$
Furthermore, we have the following estimates:
$$\|\overline h\|\leq c \sum_{i=1}^p\frac{1}{\lambda_i}\quad{\rm and }\quad\|\bar{v}\| \leq c \left(\displaystyle{
\sum_{i=1}^p (\frac{|\nabla K(a_i)|}{\lambda_i}+
\frac{1}{\lambda_i^2})+\underset{k\neq
r}{\sum}\varepsilon_{kr}\sqrt{{\rm
Log}(\varepsilon_{kr}^{-1})}}\right),$$$$
\begin{array}{rcl}J(u)&=&\displaystyle{\frac{S\sum_{i=1}^p\alpha_i^2+\alpha_0^2\|\omega\|^2}{\sqrt{S\sum_{i=1}^p\alpha_i^4K(a_i)+
\alpha_0^4\|\omega\|^2}}
\left[1-\frac{c_2\alpha_0}{\gamma_1}\sum_{i=1}^p\alpha_i\frac{\omega(a_i)}{\lambda_i}\right.}\\\\&&\displaystyle{\left.
-\frac{1}{\gamma_1}\sum_{i\neq
j}\alpha_i\alpha_jc_{ij}\varepsilon_{ij} + o\left(\sum_{i\neq
j}\varepsilon_{ij}+\sum_{i=1}^p\frac{1}{\lambda_i}\right)\right]+
\|V\|^2-\|H\|^2}
\end{array}$$
\end{prop}
{\bf Proof:} As done in \cite{Gamara prescribed} the proof is based
on performing the expansion of the functional $J$ in the set
$V(p,\varepsilon,\omega),$ to obtain  self interactions and
interactions between the bubbles, a linear form $f_{1}$ in $v$
(respectively $f_{2}$ in  $h$ ) and a  positive definite quadratic
form $Q_{1}$ in $v$ (respectively a negative definite quadratic form
$Q_{2}$ in $h$) as leading terms. Hence there is a a unique minimum
$\overline{v}$ in the space of $v's$;(respectively a unique maximum
$\overline{h}$ in the space of $h's$. Furthermore, we derive
$\|\overline{h}\|\leq c \|f_2\|$ and $\|\overline{v}\|\leq c\|
f_1\|$. The estimate of $\overline{v}$ follows from Proposition
\ref{minimum} while the estimate of $\overline{h}$ is derived from
the equivalence of the   norms  $\|\|_{\infty}$  and  $\|\|$
 in $T_{\omega}(W_u(\omega)),$ since it is a space of
finite dimension.
We also derive that $\|f_2\| = O(\sum\frac{1}{\lambda_i}),$ hence  the result follows.$\Box$\\
For the  sake of completeness of the proof one can see \cite{Bahri invariant} and \cite{Gamara prescribed}.\\
A direct consequence of the above proposition is:
\begin{cor}
Let $K$ be a $C^2$ positive function and let $\omega$  be a non
degenerate critical point of $J$ in $\Sigma^+$. Then, for each
$p\in\mathbb N^{\star}$, there is no critical points or critical
points at infinity in the set $V (p,\varepsilon,\omega)$, that means
we can construct a pseudo gradient of $J$ so that the Palais-Smale
condition is satisfied along its decreasing flow lines.
\end{cor}
The proof follows immediately from  Proposition \ref{expansionJ(w)}
and the fact that  $\omega $ is a solution of $(1.1),$ hence
strictly positive on $M.$

\section{Morse Lemma at infinity}
%%%%%%%%%%%%%%%%%%%%%%%%%%%%%%%%%%%%%%%%%%%%%%%%%%%%%%%%%%%%%%%%%%%%%%%%%%%%%%%%%
The Morse lemma at infinity establishes near the set of critical
points at infinity of the functional $J$ a change of variables in
the space $(a_{i}, \alpha_{i}, \lambda_{i}, v),$ $ 1\leq i\leq p$ to
$(\widetilde{a}_{i}, \widetilde\alpha_{i}, \widetilde\lambda_{i},
V),$  $(\widetilde\alpha_{i}=\alpha_{i})$, where $V$ is a variable
completely independent of $\widetilde{a}_{i}$ and
$\widetilde\lambda_{i}$ such that $ J(\sum \alpha_{i}
\delta_{{a}_{i},\lambda_{i}})$ behaves like $ J(\sum \alpha_{i}
\delta_{\widetilde{a}_{i}, \widetilde\lambda_{i}}) + \|V\|^{2}.$ We
 define also a pseudo-gradient for the $V$ variable
in the aim to make this variable disappear by setting
$\frac{\partial V}{\partial s}= -\nu V$ where $\nu$ is taken to be a
very large constant. Then at $s = 1,$  $V (s) = \exp(-\nu s ) V (0)$
will be as small as we wish. This shows that, in order to define our
deformation, we can work as if $V$ was zero. The deformation will be
extended immediately with the same properties to a neighborhood of
zero in the $V$ variable.

We begin by characterizing the critical points at infinity of $J$ in
the sets $V(p,\varepsilon),$ $p \geq 1$  under condition \eqref{3'}.
This characterization is obtained through the construction of a
suitable pseudogradient at infinity for the functional $J$ for which
the Palais-Smale condition is satisfied along the decreasing flow
lines as long as these flow lines do not enter in the neighborhood
of a finite number of critical points $\xi_{i};  1\leq i \leq p$ in
$\mathcal K_{2}$ or such that $(\xi_{i},... , \xi_{p})\in \mathcal
K_{1}^{+}.$
%%%%%%%%%%%%%%%%%%%%%%%%%%%%%%%%%%%%%%%%%%%%%%%%%%%%%

\subsection{ \textbf{Construction of the pseudo gradient} }
This subsection is devoted to the construction of the pseudo
gradient for the functional $J.$ It was extracted from \cite{G-H},
where a complete and detailed description of the construction of the
pseudo gradient is given.\\
 $\bullet$ In the set $V (1,\varepsilon),$ we have the following result:
\begin{prop}\label{V1}  Assume that $K$ satisfies the \textbf{$\beta-$flatness condition} and condition
$\mathbf{(C)}$   and let \\$\beta := \max\{\beta(\xi_i),\ \xi_i
\;\text{verifying}\; \eqref{3'}\}$. Then, there exists a pseudo
gradient $W$ and  a  constant $c > 0$ independent of $u =
\alpha\delta _{(a,\lambda)}\in V(1,\varepsilon),$ $\varepsilon$
small enough such that, if we denote $\overline{u}=u+ \overline{v},$
we have
\begin{enumerate}
\item $-J'(u)(W)\,\geq c (\disp\frac{|\nabla K(a)|}{\lambda}+\frac{1}{\lambda^{\beta}}).$
\item $- J'(\bar{u})(W+\disp\frac{\partial \bar{v}}{\partial(\alpha,a,\lambda)}(W))\;\geq c\,(\disp\frac{|\nabla K(a)|}{\lambda}+\frac{1}{\lambda^{\beta}}).$
\item $|W|$ is bounded. Furthermore,  $\lambda$ is an increasing function along the flow lines generated by $W$,
only if $a$ is close to a critical point $\xi_i\in \mathcal
K_1\cup\mathcal K_2.$
\end{enumerate}
\end{prop}

$\bullet$ In the set $V(p, \varepsilon),$ $p\geq 2,$ we obtain:
\begin{prop}\label{fieldW} Assume that $K$ satisfies the \textbf{$\beta-$flatness
 condition} and condition
$\mathbf{(C)}$  and let\\ $\beta := \max\{\beta(\xi_i),\ \xi_i\
\text{verifying}\  \eqref{3'}\}$.
  For any $p\geq 2$, there exists a pseudo gradient $W$ so that the following
hold:\\\vspace*{,3cm} there is a positive constant $c$ independent
of $u =\sum_{i=1}^{p}\alpha_{i}\delta _{a_{i},\lambda_{i}}\in V(p,
\varepsilon),$ $\varepsilon$ small enough such that, if we denote
$\overline{u}=u+ \overline{v},$ we have
\begin{enumerate}
 \item $-J'(u)(W)\,\geq c (\disp\sum_{i=1}^{p}\frac{|\nabla K(a_{i})|}{\lambda_{i}}+ \sum_{i=1}^{p}\frac{1}{\lambda_{i}^{\beta}}+\sum_{i\neq j}\varepsilon_{ij})$\\
 \item $- J'(\bar{u})(W+\disp\frac{\partial \bar{v}}{\partial(\alpha,a,\lambda)}(W))\;\geq c\,(\sum_{i=1}^{p}\frac{|\nabla K(a_{i})|}{\lambda_{i}}+ \sum_{i=1}^{p}\frac{1}{\lambda_{i}^{\beta}}+ \sum_{i\neq j}\varepsilon_{ij})$\\
\item $|W|$ is bounded. Furthermore, the only cases where the maximum of the $\lambda_i's$ is
not bounded is when the concentration points $(a_1, \ldots, a_p)$
satisfy: each point $a_j$ is close to a critical point $\xi_{i_j}$
of $K$ in the set $\mathcal K_1$ with $i_j \neq i_k$ for $j \neq k$
and $\varrho(\xi_{i_{1}},...,\xi_{i_{p}}) > 0$, where
$\varrho(\xi_{i_{1}},...,\xi_{i_{p}})$ is the least eigenvalue of
$M(\xi_{i_{1}},...,\xi_{i_{p}}).$
\end{enumerate}
\end{prop}

%%%%%%%%%%%%%%%%%%%%%%%%%%%%%%%%%%%%%%%%%%%%%%%%%%%%%%%%%%%%%%
\subsection{\textbf{Morse Lemma}}
Once the pseudo gradient is constructed, following \cite{Bahri
invariant} and \cite{Gamara prescribed}, we establish our Morse
Lemma at infinity: we can find a change of variables which gives the
normal form of the functional $J$ on the subsets $V(p,\varepsilon).$
We obtain the following result:
\begin{prop}\label{V2}

For $\xi \in \mathcal{K}_1\cup\mathcal K_2,$ there exists a change
of variables in the set $ \{ \alpha\delta_{(a,\lambda)} + v: a $ is
close to  $\xi \}$,  $v - \overline{v} \longmapsto V
\quad\text{and}\quad (a, \lambda) \longmapsto (\widetilde{a},
\widetilde{\lambda})$  such that in these new variables the
functional $J$  behaves as
\begin{eqnarray*}
J(\alpha\delta _{(a,\lambda)} + v) =
\frac{S}{K(\tilde{a})^{\frac{1}{2}}}
\left(1+c(1-\mu)\frac{\Gamma(\xi)}{\tilde{\lambda}^{\gamma(\xi)}}\right)+\|V\|^{2}
\end{eqnarray*}
where $\mu$ is  a small positive constant and
$$\begin{array}{rclrcl}\gamma(\xi) &=&\left\lbrace\begin{array}{cl} 2 \quad& \text{if}\;\;\xi \in \mathcal K_1\\
                                              \beta \quad&\text{if}\;\;  \xi \in \mathcal K_2
                                 \end{array}\right.; &\Gamma(\xi)&=&\begin{array}{cl} \disp-\sum_{k=1}^{2} b_{k}+\kappa' b_{0} \quad& \text{if}\;\;\xi \in \mathcal K_1\cup\mathcal K_2\\
\end{array}
\end{array}$$
\end{prop}

 The proof is similar to the one  given   in \cite{Bahri invariant,B C C H,Gamara prescribed}, so we omit it here.\\

As a consequence of Proposition \ref{V1}, we obtain:

\begin{cor} \label{p=1}Let $K$ be a positive function on $\mathbb{S}^{3}$ satisfying the \textbf{$\beta-$flatness condition} and condition
$\mathbf{(C)}.$  The only critical points at infinity in $V(1,
\varepsilon)$ are $\xi_{\infty}$ where $\xi\in\mathcal
K_1\cup\mathcal K_2$. The Morse index $i(\xi_\infty)$ of  such a
critical point is  equal to
$$\begin{array}{rcl}i(\xi_\infty) &=&  3-m(\xi)\\
\end{array}$$
\end{cor}

%%%%%%%%%%%%%%%%%%%%%%%%%%%%%%%%%%%%%%%%%%%%%%%%%%%%%%%%

If $p\geq 2,$ we have the following result:

\begin{prop}\cite{Gamara prescribed}

For any $u=\sum_{i=1}^{p}\alpha_{i}\delta _{a_{i},\lambda_{i}}\in
V(p, \varepsilon_{1}),$ $(\varepsilon_{1}< \frac{\varepsilon}{2}),$
each $a_{i}$ close to a critical point $\xi\in\mathcal K_1,$ we find
a change of variables in the space $(a_{i}, \alpha_{i}, \lambda_{i},
v),$ $ 1\leq i\leq p$ to $(\widetilde{a}_{i}, \widetilde\alpha_{i},
\widetilde\lambda_{i}, V),$  $(\widetilde\alpha_{i}=\alpha_{i})$,
such that

\begin{equation*}
\;\;\;\;\qquad J(\sum_{i=1}^{p}\alpha_{i}\delta _{a_{i},\lambda_{i}}
+ \overline{v}(\alpha, a, \lambda))= J(\sum_{i=1}^{p}
 \alpha_{i}{\delta}_{\widetilde{a}_{i}, \widetilde\lambda_{i}})
 \end{equation*}
with \begin{equation}\label{property1}  \sum_{i\neq j}
\widetilde{\varepsilon}_{ij}+
\sum_{i}\frac{1}{\widetilde\lambda_{i}^{2}}\longrightarrow 0
 \Leftrightarrow  \sum_{i\neq j}{\varepsilon}_{ij}+ \sum_{i}\frac{1}{\lambda_{i}^{2}}\longrightarrow 0.
 \end{equation}
and \begin{equation} \label{property2} \| \widetilde{a}_{i}- a_{i}
\| \longrightarrow 0\;\; \text{as}\;\; \sum_{i\neq
j}{\varepsilon}_{ij}+
\sum_{i}\frac{1}{\lambda_{i}^{2}}\longrightarrow 0
\end{equation}
\end{prop}
As a  consequence of proposition \ref{fieldW}, we obtain:

\begin{cor}\label{p} The only critical points at infinity in $ V(p, \varepsilon),$ $p\geq 2$ are:
$\xi_{\infty}=(\xi_{i_1},\ldots , \xi_{i_p} )_{\infty}$ such that
the matrix $M(\xi_{i_1},\ldots , \xi_{i_p}) $ defined in $\eqref{4}$
is positive definite, where  the  $\xi_{i_j}'s$ are critical points
of $K$ in the set $\mathcal K_1 $ and $i_j \neq i_k$ for $j \neq k.$
Such a critical point at Infinity has a Morse index equal to
$$i(\xi_{\infty})= i(\xi_{i_1},\ldots , \xi_{i_p})_{\infty}= 4p-1-
\sum_{j=1}^p m(\xi_{i_{j}})$$
\end{cor}$\Box$

\section{Proofs of Theorem \ref{existence} and Theorem \ref{multiplicity}}

\textbf{Proof of Theorem \ref{existence}}\\

Following  \cite{G-G-A} and \cite{Gamara Riahi multiplicity}, let
$\mathcal K_{\infty}$ be the set of all critical points at infinity
of $J$ and $L_0$ be their maximal Morse index given in (\ref{L0}).
We define for $0\leq l\leq L_0$ the following sets:
$$X^{\infty}_{l}=\bigcup_{\xi_\infty\in\mathcal{K}_{\infty};\;M(\xi_{\infty})\leq l}\overline {W_{u}^{\infty}(\xi_{\infty})} $$
where $W_{u}^{\infty}(\xi_{\infty})$ is the unstable manifold
associated to the critical point at infinity $\xi_{\infty}$. By a
theorem of Bahri and Rabinowitz \cite{B-R}, we have:
\begin{equation}
\overline{W_{u}^{\infty}(\xi_{\infty})}=W_{u}^{\infty}(\xi_{\infty})
\cup\bigcup_{\xi'_{\infty}<\xi_{\infty}}W_{u}^{\infty}(\xi'_{\infty})
\cup\bigcup_{\omega<\xi_{\infty}}W_{u}(\omega),
\end{equation}
where $\xi'_{\infty}$ is a critical point at infinity dominated by
$\xi_{\infty}$ and $\omega$ is a solution of $(1.1)$ dominated by
$\xi_{\infty}.$ Hence,
$$X^{\infty}_{l}=\bigcup_{\xi_{\infty}\in\mathcal K_{\infty};M(\xi_{\infty})\leq l}\Big(W_{u}^{\infty}(\xi_{\infty})
\cup\bigcup_{\omega<\xi_{\infty}}W_{u}(\omega)\Big)$$ It follows
that $X_{l}^{\infty}$ is a stratified set of top dimension $\leq l$.
Without loss of generality, we may assume it equal to $l.$ Now, we
consider  the
 cone based on  $X^{\infty}_{l}$  of vertex  $(\xi_{0})_{\infty}$ where $\xi_{0}$ is a global maximum of
$K$ on $\mathbb{S}^{3}:$
\begin{equation}
C(X^{\infty}_{l} ):= X^{\infty}_{l}\times [0, 1]\big/ (x,1)\sim
(y,1), \ x, y\in X^{\infty}_{l}
\end{equation}
 The cone $C(X^{\infty}_{l})$ is a stratified set of top
dimension $l+1$. Next,  we use the vector field $-\partial J$ to
deform $C(X^{\infty}_{l}).$ During this deformation and  based on
transversality arguments, we assume that we can avoid
 the stable manifolds of  all critical points as well as
critical points at infinity having their Morse indices greater than
$l+2$. It follows, by a Theorem of Bahri and Rabinowitz \cite{B-R},
that $C(X^{\infty}_{l})$ retracts by deformation on the set
\begin{equation}
U^{\infty}:=X^{\infty}_{l}\cup\bigcup_{M(\xi'_{\infty})=l+1}W_{u}^{\infty}(\xi'_{\infty})
\cup\bigcup_{\omega;\;\omega\ \textrm{dominated by}\
C(X^{\infty}_{l})}\;W_{u}(\omega)
\end{equation}
Now, taking $l=k-1$ and using the assumption that there are no
critical points at infinity with index k, we derive that
$C(X^{\infty}_{k-1})$ retracts by deformation onto
\begin{equation}\label{40}
Z^{\infty}_{k}:=X^{\infty}_{k-1}\cup\bigcup_{\omega;\;\omega\
\textrm{dominated by}\ C(X^{\infty}_{k-1})}\;W_{u}(\omega)
\end{equation}
Using the deformation above,  problem $\eqref{1}$ has necessary a
solution $\omega$ with $morse(\omega)\leq k$. Otherwise it follows
from $\eqref{40}$ that
$$1=\sum_{\xi_{\infty}\in\mathcal{K}_{\infty};\,M(\xi_{\infty})\leq k-1}(-1)^{M(\xi_{\infty})}=\displaystyle\sum_{\scriptsize\begin{array}{l}\xi\in \mathcal K_{2}\\ m(\xi)\geq 4-k\end{array}}(-1)^{m(\xi)+1}+\;\;
\sum_{p=1}^{l^{+}}\sum_{\scriptsize{\begin{array}{l}(\xi_{i_{1}},.,\xi_{i_{p}})\in \mathcal K^{+}\\
 i(\xi_{i_{1}},...,\xi_{i_{p}})\leq k-1 \end{array}}}
(-1)^{i(\xi_{i_{1}},...,\xi_{i_{p}})}$$   Obviously this formula
contradicts the first  assumption of the theorem.
\begin{flushright}
    $\square$
\end{flushright}
\textbf{Proof of Theorem $\ref{multiplicity}$}
 Let us denote by
$\mathcal{S}_{k}$  the set of solutions of problem $(1.1)$ having
their morse indices less than or equal to $k.$ We derive  from
$\eqref{40}$, taking the Euler characteristic of its both sides,
that:
$$1=\sum_{\xi_{\infty}\in\mathcal{K}_{\infty};\;M(\xi_{\infty})\leq k-1}(-1)^{M(\xi_{\infty})}+\sum_{\omega< C(X^{\infty}_{k-1})}(-1)^{morse(\omega)}$$
It follows that
$$\left|1-\sum_{\xi_{\infty}\in\mathcal{K}_{\infty};\,M(\xi_{\infty})\leq k-1}(-1)^{M(\xi_{\infty})}\right|\leq\#\mathcal{S}_{k}.$$
The result follows.\\

\begin{flushright}
    $\square$
\end{flushright} If we let  $k= L_{0}+ 1,$  in
Theorem \ref{existence} the second assumption of this  theorem  is
obviously satisfied and we obtain under this condition the following

\begin{cor}\label{maxexistence}
Let $K$ be  as in Theorem \ref{existence} such that: \vspace{,3cm}

$$\sum_{\xi\in \mathcal K_{2}}(-1)^{m(\xi)+1}+
\sum_{p=1}^{l^{+}}\sum_{(\xi_{i_{1}},.,\xi_{i_{p}})\in \mathcal
K_{1}^{+}}
(-1)^{i(\xi_{i_{1}},.,\xi_{i_{p}})}\neq 1$$\vspace{,3cm}\\
\vspace{,3cm}
Then, there exists at least one  solution of $\eqref{1}.$\\
\end{cor}
\begin{flushright}
    $\square$
\end{flushright}
This result generalizes the existence results due to Gamara and
Riahi in \cite{Gamara Riahi  interplay} and the multiplicity results
due to the same authors in \cite{Gamara Riahi multiplicity} and
finally recovers the existence results of Gamara and Hafassa in
\cite{G-H}.
 Moreover, if we denote $\mathcal{S}$ the set of all
the solutions of $\eqref{1},$ we obtain the following lower bound
for $\mathcal{S}$

\begin{cor} \label{maxmultiplicity}
$$\#\mathcal{S}\geq\left|1+\disp\sum_{\xi\in \mathcal
K_{2}}(-1)^{m(\xi)}-
\sum_{p=1}^{l^{+}}\sum_{(\xi_{i_{1}},.,\xi_{i_{p}})\in \mathcal
K_{1}^{+}} (-1)^{\sum_{j=1}^p m(\xi_{i_{j}})}\right|$$
\end{cor}
\begin{flushright}
    $\square$
\end{flushright}
%%%%%%%%%%%%%%%%%%%%%%%%%%%%%%%%%%%%%%%%%%%%%%%%%%%%%%%%%%%%%%%%%
\section{Appendix}
Without loss of generality, we can  assume for $p\geq 2$ that
$\lambda_{1}\leq\cdots\leq\lambda_{ p}.$ Given $N$  a large positive
constant, we define:
\begin{eqnarray}
I_1&:=&\left\{1\right\}\cup\left\{i\leq p:\lambda_k\leq
N\lambda_{k-1}\;\; \forall k\leq i\right \},\label{30}\\
I_2&:=&\left\{i\in I_1: a_i \;\hbox{{\rm is close to a critical
point}}\; \xi_{k_i}\;{\rm  satisfying}\; \eqref{3'}\;{\rm  with}\;
\beta > 2\right\}.
\end{eqnarray}
\vspace*{,3cm}
 The set $I_1$ contains the indices $i$ such that $\lambda_i$ and $\lambda_1$ are of the same
 order.\\
  We denote by $V(p, \varepsilon)_{1}$ the subset of $V(p,\varepsilon)$
composed  of the functions $u=\sum_{i=1}^{p}\alpha_{i}\delta
_{{a_{i},\lambda_{i}}}$ such that \\$\forall i\in I_1, \;\lambda_{i}
|\nabla_{\theta} K(a_{i}) | \leq 2 C^{'},
 \;\sum_{j\neq k}\varepsilon_{jk}\leq\disp\frac{C}{\lambda_{1}^{2}}, \
 C \text{ and } C' \text{ are positive constants}
 \text{ and } I_2= \emptyset.$ Following
the work done in \cite{Gamara prescribed}, we obtain the following
expansion of the functional $J$ in $ V(p, \varepsilon)_{1}.$\\

\begin{prop}\label{expansionJ}
There exists $\varepsilon_{0}>0$ such that, for any $u =
\overset{p}{\underset{j=1}{\sum}}
\alpha_{i}\delta_{a_{i,}\lambda_{i}} + v \in V(p,\varepsilon)_{1}$,
$\varepsilon < \varepsilon_{0}$, $v$ satisfying $(V_{0})$, we have
\begin{eqnarray*}
J(u) &=&\frac{\overset{p}{\underset{i=1}{\sum}}\alpha_{i}^{2}}{[\sum
\alpha_{i}^{4}K(a_{i})]^{1/2}}S \Big[1 - \frac{c}{
2S^{2}}\overset{p}{\underset{i=1}{\sum}}
\frac{\alpha_{i}^{4}}{\overset{p}{\underset{k=1}{\sum}}\alpha_{k}^{4}K(a_{k})}
\frac{\sum_{k=1}^{2} b_{k}+\kappa' b_{0}}{\lambda_{i}^{2}} \\
&& \quad + \; S^{-2}\sum_{i\neq
j}c_{0}^{4}\frac{\omega_{3}}{4}\varepsilon_{ij}
(\frac{\alpha_{i}\alpha_{j}}{\overset{p}{\underset{k=1}{\sum}}
\alpha_{k}^{2}} -
\frac{2\alpha_{i}^{3}\alpha_{j}K(a_{i})}{\overset{p}{\underset{k=1}{\sum}} \alpha_{k}^{4}K(a_{k})})\\
&&\qquad +\;\; f(v) + Q(v,v) + o(\sum_{i\neq j}\varepsilon_{ij}) +
o(\| v \|^{2}_{\theta_{1}})\Big],
\end{eqnarray*}
with
\begin{eqnarray*}
f(v) &=& - \frac{1}{\gamma_{1}}\int_{\mathbb{S}^{3}}K\big(
\overset{p}{\underset{i=1}{\sum}}
\alpha_{i}\delta_{a_{i},\lambda_{i}} \big)^{3}v\;\theta_{1}\wedge d
\theta_{1},
\end{eqnarray*}
\begin{eqnarray*}
Q(v,v) =  \frac{1}{\gamma_{2}} \|v\|^{2}_{L_{\theta_{1}}} -
\frac{3}{\gamma_{1}}
\int_{\mathbb{S}^{3}}K\overset{p}{\underset{i=1}{\sum}}\alpha_{i}^{2}
\delta_{a_{i},\lambda_{i}}^{2}v^{2}\;\theta_{1}\wedge d \theta_{1},
\end{eqnarray*}
\begin{eqnarray*}
\gamma_{1} =
S^{2}\overset{p}{\underset{i=1}{\sum}}\alpha_{i}^{4}K(a_{i})\;,\;\;\;
\gamma_{2} = S^{2}\overset{p}{\underset{i=1}{\sum}}\alpha_{i}^{2}\;.
\end{eqnarray*}

Furthermore $\| f \|_{\theta_{1}}$ is bounded

\begin{eqnarray*}
\| f \|_{\theta_{1}} &=&
O\Big(\overset{p}{\underset{i=1}{\sum}}(\frac{|\nabla
K(a_{i})|}{\lambda_{i}} + \frac{1}{\lambda_{i}^{2}}) + \sum_{i\neq
j}\varepsilon_{ij}(\log
\varepsilon_{ij}^{-1})^{\frac{1}{2}}\Big),\;\; \text{if }\; K
\;\text{satisfies }\;(\ref{3'}).
\end{eqnarray*}
\end{prop}
For a proof we refer to \cite{G-H}.
\begin{flushright}
    $\square$
\end{flushright}

%%%%%%%%%%%%%%%%%%%%%%%%%%%%%%%%%%%%%%%%%%%%%%%%%%%%%%%%%%%%%%%%%%%%%%%%%5555
Next, we will give the expansions of the gradient of the functional
$J $ which is the key of the Morse Lemma. Since the vector field $W$
is a variation of $\disp\sum_{i=1}^{p}\alpha_{i}\delta_{i}\in
V(p,\varepsilon)$ $(p\geq 2),$ hence we will expand $ J^{\prime}(u)(
\disp\lambda_j\frac{\partial\delta_j}{\partial\lambda_j}),$
$J^{\prime}(u)(
\disp\frac{1}{\lambda_j}\frac{\partial\delta_j}{\partial a_j})$ and
$ J^{\prime}(u)(\frac{1}{\lambda_{j}}(D_{j})_{k}\delta_{j}),$ for
$k=1, 2$ and $
J^{\prime}(u)(\frac{1}{\lambda_{j}^{2}}(D_{j})_{0}\delta_{j})$ in
the case where the  concentration point $a_{j},$ $j\in \{1,
2,...p\}$ is close to a critical point $\xi_j$ of $K$ verifying
$\eqref{3'}.$ We follow the lines of the method used in \cite{Gamara
prescribed} and \cite{G-G-A}. Some of the following results are
extracted from \cite{Gamara Riahi multiplicity}.\\

 For the sake of simplicity, we will use the notation  $\delta _{j}$ instead of $\delta _{a_{j},\lambda_{j}}.$
Let $u=\disp\sum_{i=1}^{p}\alpha_{i}\delta_{i}\in V(p,\varepsilon),$
we have
\begin{prop}\label{GradJ0}\cite{Gamara prescribed}
\begin{enumerate}
\item $-J^{\prime}(u)( \disp\lambda_j\frac{\partial \delta_j}{\partial\lambda_j})=
 2 J(u)\left[\disp\sum_{i\neq j}c\alpha_i\lambda_j\disp\frac{\partial\varepsilon_{ij}}{\partial \lambda_j}(1+o(1))
 - \disp\frac{\omega_3}{24}\alpha_j\disp\frac{\triangle K(a_j)}{K(a_j)\lambda_j^2}(1+o(1))\right.$\\$\left.\hspace*{4,5cm} +o\left (\sum_{i\neq j}\varepsilon_{ij}\right)\right]$
\item $-J^{\prime}(u)( \disp\frac{1}{\lambda_j}\frac{\partial \delta_j}{\partial a_j})=2 J(u)\left[ \disp\frac{\alpha_j}{K(a_j)}\frac{\omega_3}{48}\disp\frac{\nabla K(a_j)}{\lambda_j}(1+o(1))+ O\left (\sum_{i\neq j}\varepsilon_{ij}+\frac{1}{\lambda_j^2}\right)\right]$ $\Box$
\end{enumerate}
\end{prop}

 If there exists a point $a_{j},$ $j\in \{1, 2,...p\}$  close to a critical point $\xi_j$ of $K$ verifying $\eqref{3'}$,
 then the estimates in the above  proposition  can be improved see \cite{G-H} and we obtain:
\begin{prop} \label{GradJFlatness}
\begin{enumerate}
\item For $k\in\{1,2\}$
\begin{eqnarray*}
&& J^{\prime}(u)(\frac{1}{\lambda_{j}}(D_{j})_{k}\delta_{j})=
-4J(u)^{3}\alpha_{j}^{4}\frac{c_{0}^{4}}{\lambda_{j}^{\beta}}
\Big[b_{k}\int_{\mathbb{H}^{1}}\frac{|x_{k}+\lambda_{j}(a_{j})_{k}|^{\beta}}
{\Big|1+|z|^{2}-it\Big|^{6}}x_{k}\left(1+|z|^{2}\right)\theta_{0}\wedge d\theta_{0}\\
&&\hspace{2cm}+b_{0}\int_{\mathbb{H}^{1}}\frac{|t+\lambda_{j}^{2}(a_{j})_{0}+
2\lambda_{j}(x_{2}(a_{j})_{1}-x_{2}(a_{j})_{1})|^{\frac{\beta}{2}}}
{\Big|1+|z|^{2}-it\Big|^{6}}\left(x_{k}(1+|z|^{2})+(-1)^{k'}x_{k'}t\right)\theta_{0}\wedge d\theta_{0}\Big]\\
&&\hspace{2cm}+o\left(\frac{1}{\lambda_{j}^{\beta}}\right)
+O\left(\sum_{i\neq j}\varepsilon_{ij}\right)
 \end{eqnarray*}
and
\begin{eqnarray*}
&& J^{\prime}(u)(\frac{1}{\lambda_{j}^{2}}(D_{j})_{0}\delta_{j})=
-4J(u)^{3}\alpha_{j}^{4}\frac{c_{0}^{4}}{\lambda_{j}^{\beta}}
b_{0}\int_{\mathbb{H}^{1}}\frac{\Big|t+\lambda_{j}^{2}(a_{j})_{0}+
2\lambda_{j}(x_{2}(a_{j})_{1}-x_{2}(a_{j})_{1})\Big|^{\frac{\beta}{2}}}
{\Big|1+|z|^{2}-it\Big|^{6}}t\theta_{0}\wedge d\theta_{0}\\
&&\hspace{3.7cm}+o\left(\frac{1}{\lambda_{j}^{\beta}}\right)
+O\left(\sum_{i\neq j} \varepsilon_{ij}\right).
\end{eqnarray*}
\item  If we assume that $\lambda_{j}|a_{j}|\leq \mu$, where $\mu$ is a small positive constant,
then
$$ J^{\prime}(u)(\frac{1}{\lambda_{j}} \frac{\partial\delta_{j}}{\partial\lambda_{j}})=-2c_{4}J(u)\sum_{i\neq j}\alpha_{i}\lambda_{j}\frac{\partial \varepsilon _{ij}}{
\partial \lambda _{j}}+c_{5}\frac{\sum_{i=1}^{2}b_{i}+\kappa b_{0}}{\lambda_{j}^{\beta}}+
o\left(\sum_{k\neq
r}\varepsilon_{kr}+\frac{1}{\lambda_{j}^{\beta}}\right)$$$\Box$
\end{enumerate}
\end{prop}

\bibliography{mybibfile}

\end{document}